\documentclass[12pt]{amsart}
\usepackage{setspace}
\usepackage{amsfonts}
\usepackage{amsmath}
\usepackage{amssymb}
\usepackage{amsthm}
\usepackage{mathrsfs}
\usepackage[all]{xy}

\newcommand{\proofend}{\hfill \hbox{\vrule width 5pt height 5pt depth
0pt}}
\newcommand{\R}{\mathbb{R}}
\newcommand{\C}{\mathbb{C}}
\newcommand{\Z}{\mathbb{Z}}
\newcommand{\N}{\mathbb{N}}
\newcommand{\Q}{\mathbb{Q}}

\newcommand{\Oz}{\mathcal{O}}

\newcommand{\proj}{\mathbb{P}}
\newcommand{\A}{\mathbb{A}}

\newtheorem*{thm}{Theorem}

\newtheorem{lemma}{Lemma}[section]

\newtheorem{propo}[lemma]{Proposition}

\begin{document}

\title[Effectivity in Mochizuki's work on the $abc$-conjecture]{Effectivity in Mochizuki's work \\ on the $abc$-conjecture}
\author{Vesselin Dimitrov}

\maketitle

\begin{abstract}
This note outlines a constructive proof of a proposition in Mochizuki's paper {\it Arithmetic elliptic curves in general position}, making a direct use of computable non-critical Belyi maps to effectively reduce the full $abc$-conjecture to a restricted form. Such a reduction means that an effective $abc$-theorem is implied by Theorem~1.10 of Mochizuki's final IUT paper ({\it Inter-universal Teichm\"uller theory IV: log-volume computations and set-theoretic foundations}).
\end{abstract}

\section{}

Shinichi Mochizuki's proposed solution of the $abc$-conjecture revolves around the theory~\cite{mochizukietale} of the \'etale theta function of a Tate curve over a completion $F_v$ of a number field. Here of course $v$ is non-Archimedean, and moreover, the theory of~\cite{mochizukietale} has been developed under the restrictive assumption that the residue characteristic $\mathrm{char} \, k(v) \neq 2$ (cf. remark 1.10.6 (ii) in~\cite{mochizukiiv}). In the global application to the Szpiro inequality for an elliptic curve over $F$ having a good reduction outside of $S$ and a split $\mathbb{G}_m$-reduction at each place in $S$, the results of~\cite{mochizukietale} are applied at each admissible place in $S$, i.e.,  at each place in $S$ having odd residue characteristic. This\footnote{Assuming the proof is correct.} leads to an {\it \`a priori} weaker form of the full $abc$-conjecture: a Szpiro discriminant-conductor inequality for $E/F$, uniform in $F$, and restricted to elliptic curves having $|j_E|_v < C$ bounded at all  places $v$ of $F$ dividing $2$ or $\infty$. This restricted inequality is effective in its dependence on the parameter $C$.

To be a bit more precise, Theorem~1.10 of~\cite{mochizukiiv} states an explicit $2$-depleted uniform Szpiro inequality for all $E/F$ that meet certain genericity assumptions. Then the admissible choice of the auxiliary prime level $l$ and the various genericity assumptions on $E$ (big Galois image in $\mathrm{Aut} (E[l])$,\ldots) are supplied in the course of the proof of Corollary~2.2 in {\it loc.cit.} using only the assumed bounds $|j_E|_v < C$ for $v \mid \infty$, with the excluded $E$ having their heights $h(j_E)$ bounded above by an explicit function of $C$.
 This gives the $C$-restricted Szpiro bound as described above. The full strength of the $abc$-conjecture is then derived from this using the $\Sigma = V = \{2,\infty\}$ case of Theorem~2.1 of~\cite{mochizukigenell}, which for fixed values of $\varepsilon$ and $d$ yields the existence of a $C < \infty$ for which the strong $abc$-inequality for degree-$d$ points and exponent $1+\varepsilon$ is proved to be a consequence of the uniform Szpiro bound with boundedness assumptions as above: $[F:\Q] \leq d$, and  $|j_E|_v < C$ at the places $v$ of $F$ dividing $2$ or~$\infty$.

In his proof of the latter theorem Mochizuki makes an argument by contradiction, citing  compactness of $\proj^1(F_v)$. This gives the appearance of ineffectivity  of the claimed final result, Theorem~A of~\cite{mochizukiiv} (the $abc$-conjecture). The present note outlines a constructive proof, restricting for simplicity to the case
$(X,D) = (\proj^1,[0]+[1]+[\infty])$
that is actually used in the implication ``Theorem~1.10 of \cite{mochizukiiv} $\Rightarrow abc$-conjecture.'' This leads in principle to an explicit  $abc$-inequality and hence, conditionally on the correctness of Mochizuki's IUT papers, to effective Roth and Faltings theorems.

In place of the language of Arakelov theory  considered by Mochizuki we use the equivalent elementary framework of height functions attached to presented Cartier divisors as in Bombieri and Gubler~\cite{bombieri} (cf.~2.2.2 and 2.3.3 in {\it loc.cit.})  
Working in the framework of chapter~14 there we follow  its simpler notation $d(P)$ (Def.~14.3.9)  for the quantity denoted $\textrm{log-diff}_X(P)$ in Def.~1.5 (iii) in Mochizuki~\cite{mochizukigenell}, and define
$$
\mathrm{cond}_{[0]+[1]+[\infty]} := \mathrm{cond}_{[0]+[1]+[\infty]}^{\Q} : \proj^1(\bar{\Q}) \to \R^{\geq 0}
$$ 
as in~\cite{bombieri}, Example~14.4.4 (see also section~\ref{disc} below).
On the projective line $\proj^1$ over $\Z$, equip the line bundle $\Oz(1)$ with the standard metric
$$
\| ax_0 + bx_1 \| := |ax_0+bx_1| / \max(|x_0|,|x_1|)
$$
 at the Archimedean place. Denoting $\overline{\Oz(1)}$ the resulting arithmetic line bundle, the function $\mathrm{ht}_{\overline{\Oz(1)}} : \proj^1(\bar{\Q}) \to \R$ in~\cite{mochizukigenell} is just the standard absolute logarithmic Weil height denoted $h : \proj^1(\bar{\Q}) \to \R^{\geq 0}$ in Bombieri-Gubler~\cite{bombieri} (see also section~{\ref{heights}} below). Then
 the strong $abc$-conjecture of Elkies (Conjecture~14.4.12 in~\cite{bombieri}) takes the following form:
\begin{equation} \label{strongabc}
h(P) \leq (1+\varepsilon) \cdot (\mathrm{cond}_{[0]+[1]+[\infty]}(P) + d(P)) + O_{[\Q(P):\Q],\varepsilon}(1),
\end{equation}
for all $P \in \proj^1(\bar{\Q}) \setminus \{0,1,\infty\}$.

Next, for a parameter $\eta > 0$ and a finite set $\Sigma$ of places of $\Q$, consider $\mathcal{K}_{\Sigma}(\eta) \subset \proj^1(\bar{\Q})$ the set of points $P = [1:\alpha]$, where all conjugates $x =\alpha^{\sigma}$ of $\alpha \in \bar{\Q}$ fulfil
\begin{equation}  \label{cpct}
\min_{v \in \Sigma} \min \big( |x|_{v}, |1-x|_{v}, |1/x|_{v} \big) \geq \eta.
\end{equation}
Here, $|\cdot|_v$ denotes the standard absolute value on $\C_v$ (normalized by $|p|_p = 1/p$, for $p$  finite, and the ordinary absolute value on $\C$, if $v = \infty$), and the condition is independent of a choice of embeddings of $\bar{\Q}$ in $\C_v$. In Mochizuki's terminology this is a compactly bounded set supported at $\Sigma$. Our goal here is to obtain a computable function $\eta(d,\varepsilon,s)$ such that the general inequality~(\ref{strongabc}) for  $[\Q(P):\Q] \leq d$ reduces effectively to a similar inequality but with $P$ restricted to $\mathcal{K}_{\Sigma}({\eta(d,\varepsilon,|\Sigma|)})$.

\begin{thm} \label{theo}
There are computable functions
$$
\eta,c : \N \times (0,1) \times \N \to (0,1), \quad C: \N \times (0,1) \times \N \to \R
$$
such that the following is true:

Suppose $\Sigma$ a finite set of places of $\Q$ and $A : \R \times (0,1)  \to \R$  a function such that all $P \in \mathcal{K}_{\Sigma}({\eta(d,\varepsilon,|\Sigma|)})$ having $[\Q(P):\Q] \leq d$ satisfy
\begin{equation} \label{restr}
h(P) \leq (1+\varepsilon) \cdot (\mathrm{cond}_{[0]+[1]+[\infty]}(P) + d(P)) + A(d,\varepsilon),
\end{equation}
for all $(d,\varepsilon) \in \N \times (0,1)$. Then, for all $(d,\varepsilon) \in \N \times (0,1)$ and $P \in \proj^1(\bar{\Q}) \setminus \{0,1,\infty\}$ with $[\Q(P):\Q] \leq d$, it holds
\begin{eqnarray*} \label{full}
h(P) \leq (1 + \varepsilon) \cdot (\mathrm{cond}_{[0]+[1]+[\infty]}(P) + d(P)) \\  + 2A ( 150d \varepsilon^{-2} , c(d,\varepsilon,|\Sigma|) ) + C(d,\varepsilon,|\Sigma|).
\end{eqnarray*}

\end{thm}

\medskip

The meaning of ``computable'' here is straightforward. Following the outline,  using in Proposition~\ref{cond} the effective Nullstellensatz and the quantitative Chevalley-Weil theorem for curves due to Bilu, Strambi and Surroca~\cite{bilu} (see the estimation of the constant $C_1$ in the proof of Prop.~14.4.6 in~\cite{bombieri}), and estimating degrees and heights of Belyi maps on $\proj^1$ similarly to~\cite{khad}, will lead to explicit formulas for $\eta,c$, and~$C$.

 Granting this with $\Sigma = \{2,\infty\}$, an effective $abc$-statement is readily deduced from the statement of Theorem~1.10 of~\cite{mochizukiiv} by repeating the proof of Corollary~2.2 in that paper (supplying among a few other things an admissible choice of $l$ for the construction of initial $\Theta$-data), noting that the quantities ``$H_{\mathrm{unif}}$'' and ``$\eta_{\mathrm{prm}}$'' of~\cite{mochizukiiv} can be made specific, while  $\mathcal{C}_{\mathcal{K}({\eta})}$ and $\mathcal{H}_{\mathcal{K}({\eta})}$ of {\it loc.cit.} are explicit functions of $\eta$.

\section{Construction of Belyi maps} \label{belyis}

Belyi's theorem, that a regular algebraic curve $C$ over $\bar{\Q}$ can be presented as a ramified covering $f : C \to \proj^1$ with branching limited to $\{0,1,\infty\}$, is algorithmic and leads to an explicit bound on the degree and coefficients of the map $f$ in terms of  a set of defining equations for $C$. Such a bound (after applying a preliminary map to $\proj^1$) has been worked out by Khadjavi~\cite{khad}. The same applies to Mochizuki's refinement~\cite{mochizukibelyi} of Belyi's construction, or to the alternative algorithm given by Scherr and Zieve~\cite{schzieve} for the same result, upon inputting for instance the effective Riemann-Roch theorem\footnote{Giving an effective basis for the Riemann-Roch space $\mathcal{L}(D)$ in terms of the degrees and heights of defining equations of $C$ and the divisor $D$.} of Coates~\cite{coates} and Schmidt~\cite{schmidt} in the respective arguments on pages~6 of~\cite{mochizukibelyi} or~5 of~\cite{schzieve}.

Here is a direct construction avoiding Riemann-Roch in the case of the Fermat curve
$C =C_n: \, \{X^n+Y^n = Z^n\} \subset \proj^2$ over $\Q$ --- the only case used in the sequel. The Fermat curve comes with the canonical Belyi map $\pi: C_n \to \proj^1, \, [X:Y:Z] \mapsto [X^n:Z^n]$ having degree $\deg{\pi} = n^2$ and $3n$ critical points $\pi^{-1}\{0,1,\infty\}$, each of ramification index $n$; we will need to construct other Belyi maps non-critical at prescribed points, including among others the critical points of $\pi$.

To explain what is meant by effective, let $x := X/Z, y := Y/Z$ and note that every rational function  $f$ on $C_{n/\Q}$ has a unique presentation
\begin{equation} \label{presentat}
f = \frac{a_0(x)}{q_0(x)} + \frac{a_1(x)}{q_1(x)}y + \cdots + \frac{a_{n-1}(x)}{q_{n-1}(x)}y^{n-1},
\end{equation}
where $a_i \in \Z[x]$ and $q_i \in \Z[x] \setminus \{0\}$ are primitive and with no common complex root. The degrees and heights of the $x$-coordinates of the critical points of $f$ are evidently bounded as a simple function of $n$ and the degrees and heights of the polynomials $a_i,q_i$, where, as is customary, the {\it height} of a polynomial is defined to be the height of its coefficient vector as a point in a projective space.

\begin{propo} \label{belyi}
There is a computable function $B(n,d,H)$ such that the following is true. Let $S \subset \proj^1(\bar{\Q})$ be a set  with $\max_{P \in S} [\Q(P):\Q] \leq d$ and $\max_{P \in S} h(P) \leq H$. Then there exists a rational function $f : C_n \to \proj^1$ (over $\Q$) satisfying:
\begin{itemize}
\item $f$ is unramified outside of $f^{-1}\{0,1,\infty\}$;
\item $f(\pi^{-1}(S)) \not\subset \{0,1,\infty\}$;
\item the degrees and heights of the polynomials $a_i,q_i$ in the presentation~(\ref{presentat}) do not exceed $B(n,d,H)$;
\item the degree of $f$ and the heights $h(\pi(Q))$ of the points  $Q \in f^{-1}\{0,1,\infty\}$ do not exceed $B(n,d,H)$.
\end{itemize}
\end{propo}

Here,  note that upon mildly modifying the function $B$, the fourth clause is automatic from the third.

While $\pi = x^n$ ramifies at the set $\{XYZ = 0\}$, the function $x$ is unramified at $\{XZ = 0\}$ and the function $y$ is unramified at $\{YZ = 0\}$. For a local parameter along $\{Z = 0\}$  we may take $(ax-y)^{-1}$ for any rational integer $a > 1$. Indeed, multiplying the equation $x^n+y^n =1$ by $a^n$ and subtracting $y^n$, we have
$$
(ax-y)^{-1} = \frac{1+ \cdots + (ax/y)^{n-1}}{a^n/y^{n-1} - (a^n+1)y}.
$$
At each point in $\{Z = 0\}$ the denominator has a simple pole (since $a^n+1 \neq 0$), and the numerator does not vanish as $1 + a \zeta + \cdots + a^n\zeta^n \neq 0$ for $|\zeta| = 1$.

 Assume upon enlarging $S$ that $S$ is Galois stable and contains $\{0,1,\infty\}$ and, writing $S = \{0,1,\infty\} \sqcup S'$, select $a \in \Z^{> 1}$ bounded by a computable function of $n,d,H$ and subject to $y(P)/x(P) \neq a$ for all $P \in \pi^{-1}(S')$.
Then
 the map
$$
F : C_n \to \proj^1, \quad (ax(P)-y(P))^{-3 - n|S'|} \cdot x(P)y(P) \prod_{\alpha \in S'} (\pi(P) - \alpha)
$$
sends $\pi^{-1}(S)$ to $\{0\}$ and is unramified over that point. Clearly it meets condition three, and its critical value set $T \subset \proj^1(\bar{\Q})$ does not contain~$0$.

We are reduced to constructing a rational function $g : \proj^1 \to \proj^1$ (over $\Q$) such that
\begin{itemize}
\item $g$ is unramified outside of $g^{-1} \{0,1,\infty\}$;
\item $g(T) \subset \{0,1,\infty\}$;
\item $g(0) \notin \{0,1,\infty\}$;
\item the degree and height of $g$ are controlled as a function of $n,d$, and $H$.
\end{itemize}

 At this point we could follow either Lemmas 2.1, 2.1 and 2.3 in the original paper of Mochizuki~\cite{mochizukibelyi}, or Lemma~3.1, Proposition~3.2 and Proposition~3.3 of Scherr-Zieve~\cite{schzieve}. All these are manifestly algorithmic, and estimations as in Khadjavi~\cite{khad} will lead to an explicit $B(n,d,H)$.

\section{Heights} \label{heights}

In the following we consider Cartier divisors
$$
D = (U_i, f_i)_{i=1}^r \quad \textrm{on $C_{n/\Q}$},
$$
{\it with a choice of affine cover $U_1,\ldots,U_r$ of $C_n$ and local equations $f_i$}. Thus, $f_i \in K(C_n)^{\times}$, subject to $f_i/f_j \in \Oz^{\times}_{C_n}(U_i \cap U_j)$. Clearly, any Weil divisor $E$ on $C_1 \cong \proj^1$ (over $\Q$) has a distinguisthed such form via minimal equations over the standard atlas $\proj^1 \setminus \{\infty\}$ and $\proj^1 \setminus \{0\}$.

  Attached to this datum is a line bundle $O(D)$ having trivializations $O(D)|_{U_i} \cong \A^1 \times U_i$  glued by the transition function $f_i/f_j$. Its sections over $U$ are given by $r$-tuples $(h_i)_{i=1}^r$ of regular functions $h_i \in \Oz(U \cap U_i)$ on every $U \cap U_i$, linked by the patching conditions $h_i = (f_i/f_j)h_j$. The line bundle $O(D)$ comes with its tautological invertible meromorphic section $s_D = (f_i)_{i=1}^r$, which is a non-zero global section if $D$ is effective ($f_i \in \Oz(U_i)$). 

   Like in~2.2.1 of~\cite{bombieri} we consider a presentation
  $$
  \mathcal{D} = (U_1,\ldots, U_r; f_1,\ldots,f_r; D^+,\mathbf{s}^+;D^-,\mathbf{s}^-)
  $$
   to consist of a pair $D^{\pm} = (U_i, f_i^{\pm})$ of  Cartier divisors (with representatives over the same affine opens) such that $D = D^+ - D^-$ (meaning: $f_i = f_i^+ / f_i^-$) and the line bundles $O(D^{\pm})$ are base-point-free; and a choice
  \begin{eqnarray*}
  \mathbf{s}^+ = \big\{ s_w^+ =  (  h_{wi}^+ )_{i=1}^r \mid w =1,\ldots,k  \big \}, \\ \mathbf{s}^- =  \big\{ s_w^- = ( h_{wi}^- )_{i=1}^r
\mid w = 1,\ldots,l \big\}
  \end{eqnarray*}
 of generating global sections of $O(D^+)$ and $O(D^-)$.

   We say that the presentation $\mathcal{D}$ has a \emph{complexity bounded by $H \in \R$} if $n$ and the degrees and heights of all polynomials $a_i,q_i$ in the presentations~(\ref{presentat}) of the rational functions $f_j,h_{wj}^{\pm}$ are bounded by $H$.
We say that $\mathcal{D} = (D,\ldots)$ is an \emph{effective presented Cartier divisor} if $D$ is an effective Cartier divisor with a presentation $\mathcal{D}$ having $D^- = (U_i;1)_{1}^r$ (hence $O(D^-) \cong O$) and and $\mathbf{s}^- = \{1\}$. Again, an effective Weil divisor (such as $[0]+[1]+[\infty]$) on $C_1 \cong \proj^1$ has a canonical such presentation via the homogeneous minimal equation and the standard monomial basis of $\Gamma(\proj^1,O(d))$. 

   The presentation $\mathcal{D}$ gives the line bundle $O(D)$ an adelic metric by
   \begin{equation} \label{metric}
   \|s(P)\|_{\mathcal{D},v} := \min_w \max_u \Big| \frac{ss_u^-}{s_w^+}(P) \Big|_v
   \end{equation}
   for a local section $s = (h_i)_{i=1}^r$ at $P$. See~\cite{bombieri}, Prop.~2.7.11. Here $ss_u^-/s_w^+ = h_ih_{ui}^-/h_{wi}^+$ over $U_i$, and $|x|_v := |N_{K_v/\Q_p}(x)|_p^{1/[K:\Q]}$ for a place $v$ of a number field $K$ lying over the place $p$ of $\Q$, with the $p$-adic absolute value normalized by $|\cdot|_p = 1/p$ if $p$ is finite. Then the standard height on $\proj^1$ is defined by $h([\alpha:\beta]) = \sum_{v \in M_K} \max(\log{|\alpha|_v,|\beta|_v})$ over all places $M_K$ of a $K \supset \Q(\alpha,\beta)$, and the \emph{height function}  $h_{\mathcal{D}} = \mathrm{ht}_{O({\mathcal{D}})} : C_n(\bar{\Q}) \to \R$ of the metrized line bundle $O(\mathcal{D})$ is given by
   $$
   h_{\mathcal{D}}(P) = \sum_{v \in M_{\Q(P)}} -\log{\|s(P)\|_{\mathcal{D},v}}
   $$
   for any local non-vanishing section $s$ at $P$. This amounts to~2.2.2, 2.3.3 of~\cite{bombieri} while making the connection to the Arakelov theoretic framework in Mochizuki~\cite{mochizukigenell}.

\medskip

 Next, like in section~\ref{belyis}, we also say that a rational function $f$ on $C_{n/\Q}$ has a complexity bounded by $H$ if $n$ and the degrees and heights of all polynomials $a_i,q_i$ in the presentation~(\ref{presentat}) are bounded by $H$. We will  need to compare the heights $h(f(Q))$ and $h(g(Q))$ for a pair of non-constant rational functions $f$ and $g$ on $C_n$ with controlled complexities. This can be read from an algebraic dependency linking $f$ to $g$.

\begin{propo}  \label{compa}
  There are computable functions $a(H)$ and $b(H)$ such that the following is true. If $f$ and $g$ are non-constant rational functions on $C_{n/\Q}$ with complexities bounded by $H$, then
  $$
  h(f(Q)) \leq a(H)h(g(Q)) + b(H) \textrm{ for all } Q \in C_n(\bar{\Q}).
  $$ 
\end{propo}

Of course, $a(H)$ can be taken to be $\epsilon + \deg{f}/\deg{g}$ for any $\epsilon > 0$. We will not need this.

\medskip

{\it Proof. }  We look for an algebraic dependency
\begin{equation}\label{dependency}
  \sum_{i,j = 0}^{L-1} c_{ij}f^i g^j = 0
\end{equation}
with coefficients $c_{ij} \in \Z$, not all zero. Substitute the presentations~(\ref{presentat}) for $f$ and $g$, clear the $q_i(x)$ denominators, and collect the powers into monomials $x^ky^l$, $k <  2nHL \leq 2H^2L, l < n \leq H$, by using the relation $y^n = 1-x^n$.  The resulting linear system for the $L^2$ unknowns $c_{ij}$ consists of at most $2H^3L$ equations. The height of its coefficients is clearly bounded by a simple function of $H$   and so, taking $L = 4H^3$, Siegel's lemma (say) gives us a non-zero solution $(c_{ij})$ with height bounded by an explicit function of $H$. 

Since $g$ is non-constant, we have $c_{ij} \neq 0$ for some $i > 0$, and we may express~(\ref{dependency}) as
\begin{equation*}
  f^m = \sum_{ \substack{0 \leq i < m \\ -4H^3 < j < 4H^3 } } b_{ij}f^i g^j,
\end{equation*}
 for some $m < 4H^3$, with $(b_{ij})$ having rational components and an affine height bounded by an explicit function $q(H)$. Evaluating this relation at $Q$ and taking heights we get
$$
mh(f(Q)) \leq (m-1)h(f(Q)) + 4H^3 h(g(Q)) + q(H) + 6\log{H} + 5\log{2}, 
$$
giving the explicit bound with $a(H) = 4H^3$.
\proofend 

\section{Conductors and discriminants} \label{disc}

For $\mathcal{D} = (D,\ldots)$ an effective presented Cartier divisor on  $C_{n/\Q}$, the {\it conductor}
$$
\mathrm{cond}_{\mathcal{D}} : C_n(\bar{\Q}) \setminus \mathrm{supp}(D) \to \R^{\geq 0}
$$
is defined as in~\cite{bombieri}, Def.~14.4.2, with $K = \Q$ and $\lambda := -\log{\|s_D\|}_{\mathcal{D}}$ (see~(\ref{metric})):
$$
\mathrm{cond}_{\mathcal{D}}(P) := \sum_{v \in M_{\Q(P)}^{\mathrm{fin}}} \chi\big(-\log{\|s_D(P)\|_{\mathcal{D},v}} \big) \cdot \log{|1/\pi_v|_v}
$$ 
over the finite places $v$ of $\Q(P)$, where $\pi_v$ is a local parameter of $O_{\Q(P),v}$ (hence $\log{|1/\pi_v|_v} = \frac{\log{|k(v)|}}{[\Q(P):\Q]}$); and
$$
\chi(t) =  \left\{ \begin{array}{ll} 0 & \text{ if $t \leq 0$} \\ 1  & \text{ if $t > 0$.}
\end{array}    \right.
$$

 This choice of convention is made to facilitate  a direct reference to~\cite{bombieri}. Equivalently one could work with divisors on $C_{n/\Z}$ and follow~\cite{mochizukigenell}, Def.~1.5 (iii). 

  Next, as in~\cite{bombieri} Def.14.3.9, denote $d_K := \frac{1}{[K:\Q]} \log{|D_{K/\Q}|}$ the logarithmic root discriminant and, for $P \in C_n(\bar{\Q})$, write $d(P) := d_{K(P)}$. As will be crucial in the proof, it is the quantity $\mathrm{cond}_{\mathcal{D}}(P) + d(P)$ and not $\mathrm{cond}_{\mathcal{D}}(P)$ that has the good functorial property.

\begin{propo}  \label{cond}
There is a computable function $Z(H)$ such that the following is true. 

(i) Let $f : C_n \to \proj^1$ be a rational function unramified outside of $f^{-1}\{0,1,\infty\}$ and with complexity bounded by $H$. 
Then all $Q \in C_n(\bar{\Q}) \setminus f^{-1}\{0,1,\infty\}$ satisfy
\begin{eqnarray*}
d(f(Q)) + \mathrm{cond}_{[0]+[1]+[\infty]}(f(Q))  \leq d(Q) + \mathrm{cond}_{f^*([0]+[1]+[\infty])}(Q) \\ \leq d(
f(Q)) + \mathrm{cond}_{[0]+[1]+[\infty]}(f(Q)) + Z(H)
\end{eqnarray*}

(ii) Let $\mathcal{D} = (D,\ldots)$  and $\mathcal{E} = (E,\ldots)$ be effective presented Cartier divisors on $C_{n/\Q}$ with complexities bounded by~$H$ and with $D_{\mathrm{red}} = E_{\mathrm{red}}$. Then
$$
 \mathrm{cond}_{\mathcal{D}}(Q) \leq h_{\mathcal{E}}(Q) + Z(H)
$$
for all $Q \in C_n(\bar{\Q}) \setminus \mathrm{supp}(D)$. 
\end{propo}

 The first part follows from the proof of Proposition~14.4.6 of Bombieri-Gubler~\cite{bombieri} using the quantitative Chevalley-Weil theorem for curves due to Bilu, Strambi and Surroca~\cite{bilu}. (See the affine version: Theorem~1.5 of~\cite{bilu}. Use it with $\mathcal{C} = \proj^1 \setminus \{0,1,\infty\}$, $\mathcal{C}' = C_n \setminus f^{-1}\{0,1,\infty\}$, $\phi = f$ and the ``$x$'' of {\it loc.cit.} taken as, say, $j = 2^8 (x^2-x+1)^3 / x^2(1-x)^2$.) Alternatively one could follow the proof of Theorem~1.7 in Mochizuki~\cite{mochizukigenell}.
 
 The second part combines two points: $|\mathrm{cond}_{\mathcal{D}}(Q) - \mathrm{cond}_{\mathcal{E}}(Q) | \leq Z(H)/2$ and $\mathrm{cond}_{\mathcal{E}}(Q) \leq h_{\mathcal{E}}(Q) + Z(H)/2$, both
  following from Propositions~14.4.5 and~14.4.9  in Bombieri-Gubler~\cite{bombieri} and using the global section $s_E$ in computing the height $h_{\mathcal{E}}(Q)$.
   
   In both parts, the remarks  2.2.12 and 2.2.13 in~\cite{bombieri} are used, making an appeal to the effective Nullstellensatz (cf. reference~[195] in {\it loc.cit.}). 
   
   \proofend

\section{Proof of the Theorem}

Following an amplification idea of Vojta for deducing the strong $abc$-conjecture from his own conjecture with ramification for curves (see 14.4.14 and the (d) $\Rightarrow$ (a) implication in Theorem~14.4.16 of~\cite{bombieri}),  the proof is executed on the Fermat curves
$$
C_n : \, \{X^n + Y^n = Z^n\} \subset \proj^2
$$
 with their distinguished Belyi maps
$$
\pi = \pi_n : C_n \to \proj^1, \quad [X:Y:Z] \mapsto [X^n:Z^n], \quad \deg{\pi_n} = n^2.
$$
The ramification divisor of $\pi_n$ is
$$
R_{\pi_n} = (n-1) \cdot \pi_n^{-1}\{0,1,\infty\} = \Big( 1 - \frac{1}{n} \Big) \pi_n^* ([0]+[1]+[\infty]),
 $$
 giving $n \cdot R_{\pi_n} = (n-1) \pi_n^* ([0]+[1]+[\infty])$ a structure ``$n \cdot \mathcal{R}_{\pi_n}$'' as a presented Cartier divisor, of complexity bounded by a simple function of~$n$. We think of $\mathcal{R}_n$ as a presentation of $R_{\pi_n}$ as a $\Q$-divisor and define $h_{\mathcal{R}_n} := \frac{1}{n} h_{n \cdot \mathcal{R}_n}$.
Then
\begin{equation} \label{ampl}
h_{\mathcal{R}_n}(Q) = \Big( 1 - \frac{1}{n} \Big)h_{\pi^*([0]+[1]+[\infty])}(Q) = 3\Big( 1 - \frac{1}{n} \Big) h(\pi(Q)),
\end{equation}
for all $Q \in C_n(\bar{\Q})$. The essential feature here is that all critical points of $\pi_n$ lie over $\{0,1,\infty\}$ and have a high ramification index $n$. This is also the key point on page~13 (second paragraph) of Mochizuki~\cite{mochizukigenell}.

We use this as follows. Let $P \in \proj^1(\bar{\Q}) \setminus \{0,1,\infty\}$ be an arbitrary point with $[\Q(P):\Q] \leq d$, and choose $Q \in \pi_n^{-1}(P) \in C_n(\bar{\Q})$. Since $K_{C_n} \cong \pi_n^*K_{\proj^1} + R_{\pi_n}$ by Riemann-Hurwitz, we get a divisor in the class $n \cdot K_{C_n}$ with a presentation ``$n \cdot \mathcal{K}_{n}$'' of complexity bounded by a simple function of $n$, such that, defining $h_{\mathcal{K}_n} := \frac{1}{n} h_{n \cdot \mathcal{K}_n}$ again,
$h_{\mathcal{K}_{n}}(Q) = -2h(P) + h_{\mathcal{R}_{n}}(Q)$. Combined with~(\ref{ampl}) this gives
\begin{equation} \label{upst}
h_{\mathcal{K}_{n}}(Q) = \Big( 1 - \frac{3}{n} \Big) h(P).
\end{equation}

Let $\varepsilon \in (0,1)$ and choose
\begin{equation} \label{pairch}
n := 6\lceil 1/\varepsilon \rceil, \quad m := |\Sigma| \cdot dn^2 + 1,
\end{equation}
 fixing them in the following. After $m$ successive applications, Proposition~\ref{belyi} outputs $m$ morphisms $f_1 = \pi,f_2,\ldots,f_m : C_n \to \proj^1$ (over $\Q$) satisfying:
\begin{itemize}
\item[(i)] every $f_i$ is unramified outside of $f_i^{-1}\{0,1,\infty\}$;
\item[(ii)] $\pi \big( f_i^{-1}\{0,1,\infty\} \big) \cap  \pi \big( f_j^{-1}\{0,1,\infty\} \big) = \emptyset$ for $i \neq j$;
\item[(iii)] $\max_{i=1}^m \deg{f_i}$, $\max_{i=1}^m \max_{Q \in f_i^{-1}\{0,1,\infty\}} h(\pi(Q))$ and the degrees and heights of the polynomials $a_{ij},q_{ij}$ in the presentation~(\ref{presentat}) of the $f_i$ are bounded by a computable function $M(d,\varepsilon,|\Sigma|)$.
\end{itemize}

As in the introduction, fix embeddings of $\bar{\Q}$ in $\C_v$ for all $v \in \Sigma$. Consider on $\proj^1(\C_v)$ the \emph{chordal distance}
$$
\delta_v([x_0:x_1],[y_0:y_1]) := \frac{|x_0y_1-x_1y_0|_v}{\max(|x_0|_v,|x_1|_v)\max(|y_0|_v,|y_1|_v)}.
$$

 We have the Liouville bound
\begin{equation}  \label{liouville}
-\log{\delta_v(a,b)} \leq \deg{a} \cdot \deg{b} \cdot ( h(a) + h(b) + \log{2})
\end{equation}
for $a \neq b \in \proj^1(\bar{\Q})$ (cf.~Bombieri-Gubler~\cite{bombieri}, Th.~2.8.21).  From~(ii), (iii) and~(\ref{liouville}) it follows, with a computable function $1 \geq \kappa(d,\varepsilon,|\Sigma|) > 0$, that
$$
\inf_{\substack{i \neq j; \, a \in \pi (f_i^{-1}\{0,1,\infty\}) \\ b \in \pi ( f_j^{-1}\{0,1,\infty\})}} \delta_v(a, b) > 2\kappa(d,\varepsilon,|\Sigma|), \quad \textrm{for } v \in \Sigma.
$$
Consequently, for every $v \in \Sigma$, a given point of $\proj^1(\C_v)$ is within chordal distance $\kappa(d,\varepsilon,|\Sigma|)$ of at most one $\pi(f_i^{-1}\{0,1,\infty\})$. Since $m$ exceeds $|\Sigma| \cdot d \deg{\pi} \geq |\Sigma| \cdot [\Q(f_i(Q)):\Q]$,  there exists an $i$ such that the Galois orbit of $f_i(Q) \in \proj^1(\bar{\Q})$ is disjoint from $ D_{0,\kappa(d,\varepsilon)}^{(v)} \cup D_{1,\kappa(d,\varepsilon)}^{(v)} \cup D_{\infty, \kappa(d,\varepsilon)}^{(v)}$ for all $v \in \Sigma$, where $D^{(v)}_{a,r}$ is the disk $\delta_v(a,z) \leq r$ in $\proj^1(\C_v)$. Then $f_i(Q) \subset \mathcal{K}_{\Sigma}({\kappa(d,\varepsilon,|\Sigma|)})$.

Choose
$$
\epsilon := \frac{\varepsilon-\varepsilon^2}{2+8M(d,\varepsilon,|\Sigma|)^3} > 0,
$$
so that
\begin{equation} \label{later}
\frac{1+\epsilon}{1-8\epsilon \, M(d,\varepsilon,|\Sigma|)^3} < \Big( 1 - \frac{3}{n}  \Big) (1+\varepsilon)
\end{equation}
for later reference. (Recall the choice~(\ref{pairch}) making $3/n < \varepsilon/2$.) Now define the function $\eta : \N \times (0,1) \times \N \to (0,1)$ so as to have $\eta(d,\epsilon,|\Sigma|) \leq \kappa(d,\varepsilon,|\Sigma|)$ for all $\varepsilon \in (0,1)$.
Thus $f_i(Q) \in \mathcal{K}_{\Sigma}({\eta(d,\epsilon,|\Sigma|)})$.

 Since
$$
f_i^* [0] \cong f_i^* (K_{\proj^1} + [0] + [1] + [\infty]) \cong K_{C_n} + f_i^*([0]+[1]+[\infty])_{\mathrm{red}}
$$
by Riemann-Hurwitz and (i), we get a presentation ``$n \cdot \mathcal{E}_i$'' of $n \cdot f_i^*([0]+[1]+[\infty])_{\mathrm{red}}$ with complexity bounded by a computable function of $d,\varepsilon$, and $|\Sigma|$ and such that $h_{\mathcal{E}_i} := \frac{1}{n} h_{n \cdot \mathcal{E}_i}$ fulfils
\begin{equation} \label{ramh}
h(f_i(Q)) - h_{\mathcal{E}_i}(Q) = h_{\mathcal{K}_n}(Q).
\end{equation}
Then
\begin{equation}  \label{conduct}
\mathrm{cond}_{f_i^*([0]+[1]+[\infty])}  \leq h_{\mathcal{E}_i} + Z_1(d,\varepsilon,|\Sigma|)
\end{equation}
from Proposition~\ref{cond}.

\begin{quote}
\emph{Here, as in what follows, $Z_1,Z_2,\ldots$ are computable functions of their arguments, each constructed from the previous ones, starting with Propositions~\ref{compa} and~\ref{cond}.}
\end{quote}

 Applying our assumption~(\ref{restr}) at the point 
 $$
 f_i(Q) \in \mathcal{K}_{\Sigma}(\eta(d,\varepsilon,|\Sigma|))
 $$
  with $\varepsilon$ replaced by $\epsilon$, we get (noting again $[\Q(f_i(Q)):\Q] \leq  dn^2$)
\begin{eqnarray*}  \label{working}
h(f_i(Q)) \leq (1 + \epsilon) \cdot (\mathrm{cond}_{[0]+[1]+[\infty]}(f_i(Q)) + d(f_i(Q))) + A(dn^2,\epsilon) \\
\leq (1+\epsilon) \cdot ( \mathrm{cond}_{f_i^*([0]+[1]+[\infty])}(Q) + d(Q) ) + A(dn^2,\epsilon) + Z_2(d,\varepsilon,|\Sigma|) \\
\leq (1 + \epsilon) \cdot (d(Q) + h_{\mathcal{E}_i}(Q)) + A(150d/\varepsilon^2,\epsilon) + Z_3(d,\varepsilon,|\Sigma|)
\end{eqnarray*}
by Proposition~\ref{cond} and~(\ref{conduct}).

Next, by~(\ref{ramh}), (\ref{upst}) and Proposition~\ref{ampl} with $a(H) = 4H^3$, 
\begin{eqnarray*} 
 h_{\mathcal{E}_i}(Q) \leq h(f_i(Q)) \leq 4M(d,\varepsilon,|\Sigma|)^3 h(\pi(Q)) + b(M(d,\varepsilon,|\Sigma|)) \\
 \leq 8M(d,\varepsilon,\Sigma||)^3 h_{\mathcal{K}_n}(Q) + b(M(d,\varepsilon,|\Sigma|))
\end{eqnarray*}

 Using~(\ref{ramh}) and the last two displays  we conclude
\begin{eqnarray*}
h_{\mathcal{K}_n}(Q) \leq (1+\epsilon) d(Q) + \epsilon \, h_{\mathcal{E}_i}(Q)  + A(150d\varepsilon^{-2},\epsilon) + Z_3(d,\varepsilon,|\Sigma|) \\
\leq (1+\epsilon) d(Q) + 8 \epsilon M(d,\varepsilon,|\Sigma|)^3 \cdot h_{\mathcal{K}_n}(Q)+ A(150d\varepsilon^{-2},\epsilon) + Z_4(d,\varepsilon,|\Sigma|).
\end{eqnarray*}

This gives
\begin{eqnarray} \label{outs}
h_{\mathcal{K}_n}(Q) \leq \frac{1+ \epsilon}{ 1 - 8 \epsilon M(d,\varepsilon)^3} \, d(Q) + 2A(150d\varepsilon^{-2},\epsilon) \\ + Z_{5}(d,\varepsilon,|\Sigma|).
\end{eqnarray}
 Now~(\ref{upst}), (\ref{later}) and~(\ref{outs}) give
$$
h(P) \leq (1+ \varepsilon) d(Q) + 2A(150d\varepsilon^{-2},\epsilon) + Z_{6}(d,\varepsilon,|\Sigma|).
$$

 Finally, we have
\begin{eqnarray*} \label{conds}
d(Q) \leq d(Q) + \mathrm{cond}_{\pi^*([0]+[1]+[\infty])}(Q) \\
\leq d(P) + \mathrm{cond}_{[0]+[1]+[\infty]}(P) + Z_{7}(n)
\end{eqnarray*}
using Proposition~\ref{cond} again (since $\pi$ ramifies only over $[0]+[1]+[\infty]$).
 We conclude with $C(d,\varepsilon,s) = Z_{6}(d,\varepsilon,s) + 2Z_{7}(6\lceil 1/\varepsilon \rceil)$ and $c(d,\varepsilon,s) = \epsilon = (\varepsilon - \varepsilon^2) / (2 + 8M(d,\varepsilon,s)^3)$. \proofend

\end{document}